\RequirePackage{ifpdf}
\ifpdf 
\documentclass[pdftex]{sigma}
\else
\documentclass{sigma}
\fi

\begin{document}

\allowdisplaybreaks

\renewcommand{\PaperNumber}{051}

\FirstPageHeading

\ShortArticleName{Linearization of Second-Order ODEs by Generalized Sundman Transformations}

\ArticleName{Linearization of Second-Order Ordinary Dif\/ferential\\
Equations by Generalized Sundman Transformations}

\Author{Warisa NAKPIM and Sergey V. MELESHKO}

\AuthorNameForHeading{W.~Nakpim and S.V.~Meleshko}
\Address{School of Mathematics, Institute of Science, Suranaree University of
Technology,
\\ Nakhon Ratchasima, 30000, Thailand}
\Email{\href{mailto:warisa_2006@hotmail.com}{warisa{\_}2006@hotmail.com}, \href{mailto:sergey@math.sut.ac.th}{sergey@math.sut.ac.th}}

\ArticleDates{Received January 18, 2010, in f\/inal form June 03, 2010;  Published online June 15, 2010}

\Abstract{The linearization problem of a second-order ordinary dif\/ferential
equation  by the generalized Sundman transformation was considered
earlier by Duarte, Moreira and Santos  using the Laguerre form.
The results obtained in the present paper demonstrate that their
solution  of the linearization problem for a second-order ordinary
dif\/ferential equation via the generalized Sundman transformation
is not complete. We also give examples which show that the
Laguerre form is not suf\/f\/icient for the linearization problem via
the generalized Sundman transformation.}


\Classification{34A05; 34A25}

\section{Introduction}
The basic problem in the modeling of physical and other phenomena is to f\/ind
 solutions of dif\/ferential equations.
Many methods of solving dif\/ferential equations use a change of
variables that transforms a given dif\/ferential equation into another
equation with known properties. Since the class of linear equations
is considered to be the simplest class of equations, there arises the
problem of transforming a given dif\/ferential equation into a linear
equation. This problem is called a linearization
problem\footnote{The linearization problem has been studied in many
publications. A short review can be found in \cite{Meleshko,Ibragimov}.}.

The linearization problem  of a second-order ordinary dif\/ferential
equation via point transformations was solved by Sophus Lie~\cite{Lie}. He also noted that all second-order ordinary dif\/ferential equations
can be mapped into each other by means of contact transformations.
Hence, the solution of the linearization problem via contact transformations is trivial.

Comparing with the set of contact transformations the set of generalized Sundman transformations is
weaker: not any second-order ordinary dif\/ferential equation can be transformed to a~linear equation.
Hence, it is interesting to study an application of the set of genera\-li\-zed Sundman transformations to the
linearization problem of second-order ordinary dif\/ferential equations.

 The linearization problem via a generalized
Sundman transformation for second-order ordinary dif\/ferential
equations was investigated in \cite{DuarteMoreiraSantos}. The
authors of \cite{DuarteMoreiraSantos} obtained that any
second-order linearizable ordinary dif\/ferential equation which can
be mapped into the equation $u''=0$ via a generalized Sundman
transformation has to be of the form
\begin{gather}
\label{equa:MMM}
 y^{\prime \prime }+\lambda _2(x,y)y'^2+\lambda
_1(x,y)y^{\prime }+\lambda _0(x,y)=0.
\end{gather}
Using the functions
\[
\lambda_3=\lambda_{1y}-2\lambda_{2x},\qquad
\lambda_4=2\lambda_{0yy}-2\lambda_{1xy}+2\lambda_0\lambda_{2y}-\lambda_{1y}\lambda_1+2\lambda_{0y}\lambda_2+2\lambda_{2xx},
\]
they showed that equation \eqref{equa:MMM}
can be mapped into
the equation $u''=0$ via a generalized Sundman transformation if the
coef\/f\/icients $\lambda_i(x,y)$ $(i=0,1,2)$ satisfy the conditions:

(a) if $\lambda_3=0$, then $\lambda_4=0$;

(b) if $\lambda_3\neq0$, then $\lambda_4\neq0$ and the following
 equations have to be satisf\/ied
\begin{gather*}
\lambda_4^2+2\lambda_{3x}\lambda_4-2\lambda_3^2\lambda_{1x}+4\lambda_3^2\lambda_{0y}
+4\lambda_3^2\lambda_0\lambda_2-2\lambda_3\lambda_{4x}-\lambda_3^2\lambda_1^2=0,
\\
\lambda_{3y}\lambda_4+\lambda_3^2\lambda_{1y}-2\lambda_3^2\lambda_{2x}-\lambda_3\lambda_{4y}=0.
\end{gather*}

The generalized Sundman transformation was also applied
\cite{EulerWolfLeachEuler,NakpimMeleshko} for obtaining necessary and
suf\/f\/icient conditions for a third-order ordinary dif\/ferential
equation to be equivalent to a linear equation in the Laguerre form.
Some applications of the generalized Sundman transformation to
ordinary dif\/ferential equations were considered in~\cite{Berkovich2001} and earlier papers, which are summarized
in the book~\cite{Berkovich2002}.

According to the Laguerre theorem in any linear ordinary dif\/ferential equation the two terms of order
below next to highest can be simultaneously removed by a point transformation. For example, the Laguerre
form of a second-order ordinary dif\/ferential equation is the
linear equation $u''=0$.
For obtaining
this form, several point transformations are applied
consecutively. Since the composition of point transformations is a
point transformation, the f\/inal transformation is again a point
transformation. This is not the case for generalized Sundman
transformations: the composition of a point transformation and a
generalized Sundman transformation is not necessarily a
generalized Sundman transformation. Hence, for the linearization
problem via generalized Sundman transformations it is not
suf\/f\/icient to use the Laguerre form.

In this paper, we demonstrate that the solution of the
linearization problem via the gene\-ra\-li\-zed Sundman transformation
of second-order ordinary dif\/ferential equations given in
\cite{DuarteMoreiraSantos} only gives particular criteria
for linearizable equations.
Complete analysis of the compatibility of arising equations is given for the case $F_x=0$.

\section{Generalized Sundman transformations}

A generalized Sundman transformation is a non-point
transformation def\/ined by the formulae
\begin{gather}
\label{equa:01}
u(t) = F(x,y),\qquad dt = G(x,y)dx,\qquad F_yG\ne 0.
\end{gather}
Let us explain how the generalized Sundman transformation maps one
function into another.

Assume that $y_0(x)$ is a given function. Integrating
the second equation of (\ref{equa:01}),
one obtains $t=Q(x)$, where
\[
Q(x)= t_0+\int_{x_0}^{x}  G(s,y_0(s))\,ds
\]
with some initial conditions $t_0$ and $x_0$.
Using the inverse function
theorem, one f\/inds $x=Q^{-1}(t)$. Substituting $x$ into the
function $F(x,y_0(x)) $, one gets the transformed function
\[
u_0(t) = F\left(Q^{ - 1} (t),y_0 (Q^{ - 1} (t))\right).
\]

Conversely, let $u_0(t)$ be a given function of $t$. Using the
inverse function theorem, one solves the equation
\[
u_0 (t) = F(x,y)
\]
with respect to $y$: $y=\phi \left( {x,t}\right)$. Solving the
ordinary dif\/ferential equation
\[
\frac{dt} {dx} = G(x,\phi (x,t)),
\]
one f\/inds $t=H(x)$. The function $H(x)$ can be written as an action of a functional $H={\cal L}(u_0)$.
Substituting $t=H(x)$ into the function $\phi(x,t)$, the transformed function $y_0(x)=\phi
(x,H(x))$ is obtained.

Notice that for the case
$G_y=0$ the action of the functional ${\cal L}$ does not depend on the function $u_0(t)$.
In this case the generalized Sundman transformation becomes a point transformation.
Conversely, since for a point transformation the value $dt$ in the generalized Sundman transformation
is the total dif\/ferential of $t$, then the compatibility condition for $dt$ to be a~total dif\/ferential
leads to the equation $G_y=0$. Hence, the generalized Sundman transformation is a~point transformation if and only if $G_y=0$.

Formulae (\ref{equa:01}) also allow  us to obtain the derivatives
of $u_0(t)$ through the derivatives of the function $y_0(x)$, and
vice versa.

Hence, using transformation (\ref{equa:01}), one can relate the
solutions of two dif\/ferential equations
$Q(x,y,y',\ldots,y^{(n)})=0$ and $P(t,u,u',\ldots,u^{(n)})=0$.
Therefore the knowledge of the general solution of one of them
gives the general solution of the other equation, up to solving
one ordinary dif\/ferential equation of f\/irst-order and f\/inding two
inverse functions.

\section{Necessary conditions}

We start with obtaining necessary conditions for the linearization
problem.

First, one f\/inds the general form of a second-order ordinary
dif\/ferential equation
\[
y''  = H\big( x,y,y' \big),
\]
which can be mapped via a generalized Sundman transformation into
the linear equation
\begin{gather}
\label{equa:KL}
u''+\beta u'+\alpha u=\gamma,
\end{gather}
where $\alpha(t)$, $\beta(t)$ and $\gamma(t)$ are some functions.
Notice that the Laguerre form of a linear second-order ordinary
dif\/ferential equation corresponds to $\alpha=0$, $\beta=0$ and
$\gamma=0$.

The function $u$ and its derivatives $u'$ and $u''$ are def\/ined by the f\/irst formula (\ref{equa:01})
and its derivatives with respect to $x$:
\begin{gather}
u^{\prime }G=F_x+F_yy^{\prime },\nonumber
\\
u^{\prime \prime }G^2+u^{\prime }(G_x+G_yy^{\prime })=F_yy^{\prime
\prime }+2F_{xy}y^{\prime }+F_{yy}y^{\prime 2}+F_{xx}.\label{equa:02}
\end{gather}
The independent variable $t$ is def\/ined by the functional  ${\cal L}(u)$.
As noted above, if $G_y\neq 0$, then the action of the functional ${\cal L}$ depends on the function $u$.
Hence, if one of the coef\/f\/icients (\ref{equa:KL})
is not constant and $G_y\neq 0$, then the substitution of $t$ into equation
(\ref{equa:KL})
gives a functional equation. Since the case $G_y=0$ reduces the generalized Sundman
transformation to a point
transformation\footnote{Later it will be shown that in our study for $F_x=0$
this case is automatically excluded from
consideration because in the process of studying
the compatibility this case leads to the conditions $\lambda_3=0$ and $\lambda_4=0$,
which were considered in~\cite{DuarteMoreiraSantos}.},
the irreducible generalized Sundman transformation maps equation~(\ref{equa:KL})
into a~dif\/ferential equation only for constant coef\/f\/icients
$\alpha$, $\beta$ and $\gamma$.
Thus, f\/inding the derivatives~$u^{\prime }$,~$u^{\prime \prime }$ from
(\ref{equa:02}), and substituting them into (\ref{equa:KL}) with constant coef\/f\/icients, one
has the following equation
\begin{gather}
\label{equa:03}
y^{\prime \prime  }+\lambda _2(x,y)y'^2+\lambda _1(x,y)y^{\prime
}+\lambda _0(x,y)=0,
\end{gather}
where the coef\/f\/icients  $\lambda _i(x,y)$ $(i=0,1,2)$ are related to the functions
$F$ and $G$:
\begin{gather}
 \label{equa:04} \lambda _2 = (F_{yy}G - F_yG_y)/K, \\
   \label{equa:05}\lambda _1= (2F_{xy}G- F_xG_y -F_yG_x + F_y\beta G^2)/K, \\
  \label{equa:06} \lambda_0=(F_{xx}G - F_xG_x + F_x\beta G^2 + \alpha FG^3- G^3\gamma)/K,
\end{gather}
where $K=GF_y\neq0$.

Equation (\ref{equa:03}) presents the necessary form of a
second-order ordinary dif\/ferential equation which can be mapped
into a linear equation (\ref{equa:KL}) via a generalized Sundman
transformation.

\section{Suf\/f\/icient conditions}

 For obtaining suf\/f\/icient conditions, one has to solve the
compatibility problem (\ref{equa:04})--(\ref{equa:06}), con\-si\-de\-ring
(\ref{equa:04})--(\ref{equa:06}) as an overdetermined system of
partial dif\/ferential equations  for the functions $F$ and $G$ with
the given coef\/f\/icients $\lambda _i(x,y)$ $(i=0,1,2)$. Notice that
the compatibility conditions (\ref{equa:04})--(\ref{equa:06}) for
the particular case $\alpha=0$, $\beta=0$ and $\gamma=0$ were
obtained in~\cite{DuarteMoreiraSantos}. This case corresponds to the Laguerre form of a
linear second-order ordinary dif\/ferential equation. It is shown
here that for the linearization problem via generalized Sundman
transformations it is not suf\/f\/icient to use the Laguerre form.

The compatibility analysis depends on the value of $F_x$. A
complete study of all cases is cumbersome. Here a
complete solution is given for the case where $F_x=0$.

Solving  equations (\ref{equa:04})--(\ref{equa:06}) with respect to
$F_{yy}$, $\beta$ and $\gamma$, one f\/inds
\begin{gather}
\label{equa:07}
F_{yy}= (G_yF_y +F_yG\lambda_2)/G,
\\
\label{equa:08}
\beta= (  G_x  + G\lambda_1)/ G^2,
\\
\label{equa:09}
\gamma=(- F_y\lambda_0 + \alpha FG^2)/G^2 .
\end{gather}
Since $F_x=0$, then dif\/ferentiating $F_{yy}$ with respect to $x$, one obtains
\begin{gather}
\label{equa:100}
GG_{xy}-G_xG_y + \lambda_{2x}G^2=0.
\end{gather}
Dif\/ferentiating (\ref{equa:08}) and (\ref{equa:09}) with respect
to $x$ and $y$, one obtains the following equations
\begin{gather}
\label{equa:097}
 G_{xx}= (2G_x^2 + G_xG\lambda_1 -
\lambda_{1x}G^2)/G,
\\
\label{equa:096}
 G_{xy}= G\lambda_3- G_y\lambda_1 ,
\\
\label{equa:00} 2G_x\lambda_0 - \lambda_{0x}G=0,
\\
\label{equa:MN} \alpha= ( - G_y\lambda_0 + G(\lambda_{0y} +
\lambda_0\lambda_2))/G^3,
\end{gather}
where
\[
\lambda_3=\lambda_{1y} - 2\lambda_{2x}.
\]
Substituting (\ref{equa:096}) into (\ref{equa:100}), this becomes
\begin{gather}
\label{equa:098}   G_xG_y + G_yG\lambda_1 - G^2(\lambda_{2x} +
\lambda_3)=0.
\end{gather}
Comparing the mixed derivatives $(G_{xy})_x=(G_{xx})_y$, one
obtains the equation
\begin{gather}
\label{equa:000}
G_x\lambda_3 - G(\lambda_{2xx} + \lambda_{2x}\lambda_1 +
\lambda_{3x})=0.
\end{gather}
Dif\/ferentiating $\alpha$ with respect to $x$ and $y$, one has
\begin{gather}
\label{equa:001}
 2G_x(\lambda_{0y} + \lambda_0\lambda_2) +G_y( \lambda_{0x} +
2\lambda_0\lambda_1) - G(\lambda_{0xy} +\lambda_{0x}\lambda_2
+4\lambda_{2x}\lambda_0 + 2\lambda_0\lambda_3)=0 ,
\\
\label{equa:002}
2GG_{yy}\lambda_0-6G_y^2\lambda_0 + 2G_yG( 3\lambda_{0y} +
2\lambda_0\lambda_2) - G^2(\lambda_4
+ 2\lambda_5-\lambda_1\lambda_3)=0,
\end{gather}
where
\begin{gather*}
\lambda_4=2\lambda_{0yy}-2\lambda_{1xy}+2\lambda_0\lambda_{2y}-\lambda_{1y}\lambda_1+2\lambda_{0y}\lambda_2+2\lambda_{2xx},
\\
\lambda_5=\lambda_{2xx} + \lambda_{2x}\lambda_1 + \lambda_{3x} +
\lambda_1\lambda_3.
\end{gather*}

Further analysis of the compatibility depends on $\lambda_3$.

\subsection[Case lambda_3 not= 0]{Case $\boldsymbol{\lambda_3 \neq 0}$}

From equations (\ref{equa:000}), one f\/inds
\begin{gather}
\label{equa:AAA}
G_x= G(\lambda_{2xx} + \lambda_{2x}\lambda_1 +
\lambda_{3x})/\lambda_3.
\end{gather}
Substituting $G_x$ into equations (\ref{equa:00}),
(\ref{equa:098}),  (\ref{equa:097}) and (\ref{equa:096}), one
obtains the equations
\begin{gather}
\label{equa:1b}
 \lambda_{0x}= 2\lambda_0( - \lambda_1\lambda_3 +
\lambda_5)/\lambda_3,
\\
\label{equa:1a}
\lambda_{2xxy}= - \lambda_{2xy}\lambda_1- \lambda_{3xy} -
2\lambda_{2x}^2 - 2\lambda_{2x}\lambda_3
 - \lambda_{3y}\lambda_1 + (\lambda_{3y}\lambda_5)\lambda_3^{-1},
\\
\label{equa:1d}
\lambda_{2xxx}=- \lambda_{3xx}-
\lambda_{1x}\lambda_{2x} - \lambda_{1x}\lambda_3 +
\lambda_{2x}\lambda_1^2 + \lambda_1^2\lambda_3
-2\lambda_1\lambda_5+\lambda_3^{-1}\lambda_5(\lambda_{3x}+\lambda_5) ,
\\
\label{equa:1c}
 G_y\lambda_5-
G\lambda_3(\lambda_{2x}+\lambda_3)=0.
\end{gather}

\subsubsection[Case lambda_5 not= 0]{Case $\boldsymbol{\lambda_5 \neq 0}$}

Equation (\ref{equa:1c}) gives
\begin{gather}
\label{equa:ZZZ} G_y=
G\lambda_3(\lambda_{2x}+\lambda_3)/\lambda_5.
\end{gather}
Substituting $G_y$ into equations (\ref{equa:096}),
(\ref{equa:001}) and (\ref{equa:002}) and comparing the mixed
derivatives $(G_{x})_y=(G_{y})_{x}$, one gets
\begin{gather}
\lambda_3\lambda_5(6\lambda_{0y}\lambda_{2x}+2\lambda_{2xy}\lambda_0
+4\lambda_{2x}\lambda_0\lambda_2+2\lambda_{3y} \lambda_0+4\lambda_0\lambda_2\lambda_3
+\lambda_1\lambda_5)\nonumber\\
\qquad{}
-\lambda_3^2( 6\lambda_{2x}^2\lambda_0+
12\lambda_{2x}\lambda_0\lambda_3- 6\lambda_{0y}\lambda_{5}+
6\lambda_0\lambda_3^2)- \lambda_4\lambda_5^2 - 2\lambda_5^3=0.\label{equa:1g}
\end{gather}

\subsubsection[Case lambda_5=0]{Case $\boldsymbol{\lambda_5 = 0}$}

Equations (\ref{equa:1b}), (\ref{equa:1c}), (\ref{equa:1a}),
(\ref{equa:1d}), (\ref{equa:001}) and (\ref{equa:002}) become
\begin{gather}
\label{equa:1h} \lambda_{0x}=  - 2\lambda_0\lambda_1,
\\
\label{equa:1hhh} \lambda_{2x}=-\lambda_3,
\\
\label{equa:1i}
 2GG_{yy}\lambda_0-6G_y^2\lambda_0 + 2G_yG( 3\lambda_{0y} +
2\lambda_0\lambda_2) - G^2(\lambda_4 - \lambda_1\lambda_3)=0.
\end{gather}

If $\lambda_0 \neq0$, then equation (\ref{equa:1i}) def\/ines
\begin{gather}
\label{equa:ABC} G_{yy}=(6G_y^2\lambda_0 - 2G_yG(  3\lambda_{0y}
+2\lambda_0\lambda_2 ) + G^2( \lambda_4 - \lambda_1\lambda_3
))/(2G\lambda_0).
\end{gather}
In this case, $(G_{yy})_x=(G_{xy})_y$ and
$(G_{x})_{yy}=(G_{yy})_x$ are satisf\/ied. Hence, there are no other
compatibility conditions. Thus, if $\lambda_3\neq 0$ , $\lambda _5=0$ and
$\lambda _0\neq0$, then conditions (\ref{equa:1h}) and
(\ref{equa:1hhh}) are suf\/f\/icient for equation (\ref{equa:03}) to
be linearizable by a generalized Sundman transformation.

If $\lambda_0 =0$, there is no other conditions.

\begin{remark}
If  $\lambda_5 = 0$, equations (\ref{equa:1b}),
 (\ref{equa:1a}),
(\ref{equa:1d}), (\ref{equa:1c})
and (\ref{equa:1g}) become conditions (\ref{equa:1h}) and
(\ref{equa:1hhh}) respectively.
\end{remark}

Thus, suf\/f\/icient conditions for equation (\ref{equa:03}) in the case $\lambda_3\neq 0$
to be linearizable by  generalized Sundman transformations are (\ref{equa:1b}), (\ref{equa:1a}), (\ref{equa:1d})
and (\ref{equa:1g}).

\subsection[Case lambda_3=0]{Case $\boldsymbol{\lambda_3 = 0}$}

Notice that the particular case $\lambda_3=0$ and $\lambda_4=0$
was studied in \cite{DuarteMoreiraSantos}. Here the case $\lambda_3=0$ and
$\lambda_4\neq0$ is considered.

Equation (\ref{equa:002}) for $\lambda_3=0$ becomes
\begin{gather}
\label{equa:2d}
2GG_{yy}\lambda_0-6G_y^2\lambda_0 + 2G_yG( 3\lambda_{0y} +
2\lambda_0\lambda_2) - G^2\lambda_4 =0.
\end{gather}

The assumption $\lambda_0=0$ leads to the contradiction that
$\lambda_4=0$. Hence, one has to assume that $\lambda_0\neq0$.

Equations (\ref{equa:000}), (\ref{equa:00}) and (\ref{equa:001})
become
\begin{gather}
\label{equa:2a}
 \lambda_{2xx}= -\lambda_{2x}\lambda_1,
\\
\label{equa:2b} G_x= (G\lambda_{0x})/(2\lambda_0),
\\
\label{equa:2c} G_y\lambda_0\lambda_6 - G(\lambda_{6y}\lambda_0-
\lambda_{0y}\lambda_6 )=0,
\end{gather}
where
\[
\lambda_6=\lambda_{0x} + 2\lambda_0\lambda_1.
\]
Substituting $G_x$ into equations (\ref{equa:096})  and
(\ref{equa:097}), one gets
\begin{gather}
\label{equa:2e}
 \lambda_{6y}= (\lambda_{0y}\lambda_6 +
2\lambda_{2x}\lambda_0^2)/\lambda_0,
\\
\label{equa:2f}\lambda_{6x}=(3\lambda_6( \lambda_6-
2\lambda_0\lambda_1 ))/(2\lambda_0).
\end{gather}

\subsubsection[Case lambda_6 not= 0]{Case $\boldsymbol{\lambda_6 \neq0}$}

From equations (\ref{equa:2c}), one f\/inds
\[
G_y=G( - \lambda_{0y}\lambda_6 +
\lambda_{6y}\lambda_0)/(\lambda_0\lambda_6).
\]
Substituting $G_y$ into equations  (\ref{equa:096}) and
(\ref{equa:2d}),  and comparing the mixed derivatives
$(G_{x})_y=(G_{y})_{x}$, one obtains
\begin{gather}
\label{equa:2g}
 \lambda_{4x}= ( - 24\lambda_{2x}^2\lambda_0^3 -
4\lambda_0\lambda_1\lambda_4\lambda_6 +
\lambda_4\lambda_6^2)/(2\lambda_0\lambda_6).
\end{gather}

\subsubsection[Case lambda_6 =0]{Case $\boldsymbol{\lambda_6 =0}$}

In this case equation (\ref{equa:2c}) is satisf\/ied. One needs to
check the only condition $(G_{yy})_x=(G_{x})_{yy}$, which is
\begin{gather}
\label{equa:2hh} \lambda_{4x}=  - 2\lambda_1\lambda_4.
\end{gather}
Equation (\ref{equa:2e}) becomes
\begin{gather}
\label{equa:2h} \lambda_{2x}=0.
\end{gather}

\begin{remark}
If  $\lambda_6 = 0$, equation (\ref{equa:2e}) becomes a condition (\ref{equa:2h}).
\end{remark}

All obtained results can be summarized in the theorem.

\begin{theorem}\label{theorem}
Sufficient conditions for equation \eqref{equa:03}
to be linearizable via a generalized Sundman transformation
with $F_x=0$ are as follows.
\begin{enumerate}\itemsep=0pt
\item[$(a)$] If $\lambda_3\neq 0$, then the
conditions are \eqref{equa:1b}, \eqref{equa:1a}, \eqref{equa:1d}
and \eqref{equa:1g}.

\item[$(b)$] If $\lambda_3= 0$, $\lambda_6\neq0$, then the conditions are
\eqref{equa:2a}, \eqref{equa:2e}, \eqref{equa:2f} and
\eqref{equa:2g}.

\item[$(c)$]  If $\lambda_3= 0$, $\lambda_6=0$, then the conditions are
\eqref{equa:2a}, \eqref{equa:2e}, \eqref{equa:2f} and
\eqref{equa:2hh}.
\end{enumerate}
\end{theorem}

\begin{remark}
These conditions  extend the criteria obtained in
\cite{DuarteMoreiraSantos} to the case $\alpha, \beta,\gamma\neq 0$ in (3), for restricted ($F_x=0$) generalized
Sundman transformations.
\end{remark}

\begin{remark}
Notice that a discussion of the case
$\lambda_i=\lambda_i(y)$ $(i=0,1,2)$ and $G_x=0$, is also given in~\cite{Berkovich2001}.
\end{remark}

\begin{remark}
Recall S.~Lie's results \cite{Lie} on linearization of a second-order
ordinary dif\/ferential equation via a change of the independent and
dependent variables (point transformations).
The necessary form of a linearizable equation $y''=f(x,y,y')$ has to be the following
 form
\begin{gather}
\label{eq:Oct1604.4}
y^{\prime \prime }+a(x,y)y^{\prime }{}^3+b(x,y)y^{\prime
}{}^2+c(x,y)y^{\prime }+d(x,y)=0.
\end{gather}
Equation (\ref{eq:Oct1604.4}) is linearizable if and only if its coef\/f\/icients
satisfy the conditions
\begin{gather}
3a_{xx}-2b_{xy}+c_{yy}-3a_xc+3a_yd+2b_xb-3c_xa-c_yb+6d_ya=0, \nonumber\\
b_{xx}-2c_{xy}+3d_{yy}-6a_xd+b_xc+3b_yd-2c_yc-3d_xa+3d_yb=0.\label{eq:Oct1604.11}
\end{gather}

Despite that the form (\ref{equa:03}) is a particular case of (\ref{eq:Oct1604.4}),
suf\/f\/icient conditions of linearization via point transformations (\ref{eq:Oct1604.11}) and the generalized Sundman transformation dif\/fer.
Hence, the second class of equations is not contained in the f\/irst class due
to dif\/ferences in conditions on arbitrary elements of the classes. At
the same time, these classes have a nonempty intersection.
\end{remark}

\section{Examples}

\begin{example}
Consider the nonlinear ordinary dif\/ferential
equation
 \begin{gather}
 \label{ex:01}
y''+(1/y)y'^2 + yy'+1/2=0.
 \end{gather}
Since this equation does not satisfy Lie criteria for linearization \cite{Lie}
it is not linearizable by point transformations.
Equation (\ref{ex:01}) is of the form (\ref{equa:03}) with
coef\/f\/icients
\begin{gather}
 \lambda_2 = 1/y, \qquad  \lambda_1 = y,  \qquad \lambda_0 =1/2.
 \label{ex:BB}
\end{gather}
One can check that the coef\/f\/icients (\ref{ex:BB}) obey the
conditions (\ref{equa:1b}), (\ref{equa:1a}), (\ref{equa:1d}) and
(\ref{equa:1g}). Thus, equation (\ref{ex:01}) is linearizable via
generalized Sundman transformation.

For f\/inding the functions $F$ and $G$ one has to solve equations
(\ref{equa:07}), (\ref{equa:AAA})  and (\ref{equa:ZZZ}), which
become
\begin{gather*}
   F_{x} =0, \qquad F_{yy}=(2F_y)/y,  \qquad G_x =0, \qquad G_y=
G/y.
\end{gather*}

We take the simplest solution,  $F=y^3$ and $G=y$, which satisf\/ies
(\ref{equa:07}), (\ref{equa:AAA})  and (\ref{equa:ZZZ}). One
obtains the transformation
\begin{gather}
\label{ex:BBN} u = y^3,\qquad dt = ydx.
\end{gather}
Equations (\ref{equa:08}), (\ref{equa:09}) and (\ref{equa:MN})
give
\begin{gather*}
 \beta =1, \qquad \gamma=-3/2, \qquad \alpha=0.
\end{gather*}
Hence equation (\ref{ex:01}) is mapped by the transformation
(\ref{ex:BBN}) into the linear equation
 \begin{gather}
 \label{ex:16}
u''+u'+3/2=0.
 \end{gather}
The general solution of equation (\ref{ex:16}) is
\[
u=c_1+c_2e^{-t}-3t/2,
\]
where $c_1$, $c_2$  are arbitrary constants. Applying the generalized Sundman
transformation~(\ref{ex:BBN}) to equation~(\ref{ex:01}) one
obtains that the general solution of equation (\ref{ex:01}) is
\[
y(x)=(c_1+c_2e^{-\phi(x)} -3\phi(x)/2)^{1/3},
\]
where the function $t=\phi(x)$ is a solution of the equation
\[
 \frac{dt}{dx} =\big(c_1+c_2e^{-t} -3t/2\big)^{1/3}.
\]
For example, if $c_1=c_2=0$, then one obtains the solution of
equation (\ref{ex:01}):
\[
y=(-x)^{1/2}.
\]
\end{example}

\begin{example}
Consider the nonlinear ordinary dif\/ferential
equation
  \begin{gather}
 \label{ex:02}
y''+xy'^2 + yy'+1/{e^{2xy}}=0.
 \end{gather}
Equation (\ref{ex:02}) is of the form~(\ref{equa:03}) with the
coef\/f\/icients
\begin{gather}
 \lambda_2 = x, \qquad  \lambda_1 = y, \qquad \lambda_0 = 1/{e^{2xy}}.
 \label{ex:CC}
\end{gather}
One can check that the coef\/f\/icients (\ref{ex:CC}) do not satisfy
the conditions of linearizability by point transformations, but
they obey the conditions (\ref{equa:1h}) and
(\ref{equa:1hhh}). Thus, equation~(\ref{ex:02}) is linearizable
via a generalized Sundman transformation.

For f\/inding the functions $F$ and $G$ one has to solve equations
(\ref{equa:07}), (\ref{equa:AAA})  and (\ref{equa:ABC}), which
become
\begin{gather*}
    F_{x} =0, \qquad F_{yy}=(G_yF_y + F_yGx)/G, \\
     G_x =-yG, \qquad G_{yy}= (3G_y^2 + 4G_yGx + 2G^2x^2)/G.
\end{gather*}

We take the simplest solution,  $F=y$ and $G=e^{-xy}$, which
satisf\/ies (\ref{equa:07}), (\ref{equa:AAA})  and (\ref{equa:ABC}).
The linearizing generalized Sundman transformation is
\begin{gather}
\label{ex:15} u = y,\qquad dt = e^{-xy}dx.
\end{gather}
Equations (\ref{equa:08}), (\ref{equa:09}) and (\ref{equa:MN})
give
\begin{gather*}
   \beta =0, \qquad \gamma=-1, \qquad \alpha=0.
\end{gather*}
Hence equation (\ref{ex:02}) is mapped by the transformation
(\ref{ex:15}) into the linear equation
 \begin{gather}
 \label{ex:BBM}
u''+1 =0.
 \end{gather}
The general solution of equation (\ref{ex:BBM}) is
\[
u=-t^2/2+c_1t+c_2,
\]
where $c_1$, $c_2$  are arbitrary constants. Applying the generalized Sundman
transformation (\ref{ex:15}) to equation~(\ref{ex:02}) one obtains
that the general solution of equation~(\ref{ex:02}) is
\[
y(x)=-\phi(x)^2/2+c_1\phi(x)+c_2,
\]
where the function $t=\phi(x)$ is a solution of the equation
\[
 \frac{dt}{dx} =e^{-x(-t^2/2+c_1t+c_2)}.
\]
\end{example}

\begin{example}
Consider the nonlinear second-order ordinary dif\/ferential
equation
\begin{gather}
\label{ex:00}y^{\prime \prime }+\mu _3y^{k_3}y'^2+\mu _2y^{k_2}y^{\prime }+\mu
_1y^{k_1}=0,
\end{gather}
where $k_1$, $k_2$, $k_3$, $\mu _1$, $\mu _2$ and $\mu_3\neq 0$ are arbitrary constants. The
Lie criteria \cite{Lie} show that the nonlinear equation
(\ref {ex:00}) is linearizable by a point transformation if and
only if $\mu_1=0$ and $\mu _2=0$.

From equation (\ref {ex:00}), the coef\/f\/icients are
\begin{gather*}
\lambda _0=\mu _1y^{k_1},\qquad \lambda _1=\mu_2y^{k_2},\qquad \lambda
_2=\mu_3y^{k_3},\qquad \lambda
_3=\mu_2k_2y^{k_2}/y, \nonumber\\
\lambda _4=2\mu _1y^{(k_1+k_3)+1}(k_1\mu_3+k_3\mu_1)+2\mu_1y^{k_1}(k_1^2-k_1)-k_2\mu_2^2y^{2k_2+1}/y^2,\nonumber\\ \lambda _5=k_2\mu_2^2y^{2k_2}/y.
\end{gather*}

If $\mu_2\neq0$ and $\mu_1=0$, then $\lambda _3\neq0$ and
$\lambda _5\neq0$. One can check that the coef\/f\/icients obey the
conditions (\ref{equa:1b}), (\ref{equa:1a}), (\ref{equa:1d}) and
(\ref{equa:1g}). Thus, equation
\begin{gather}
\label{ex:VV}y^{\prime \prime }+\mu _3y^{k_3}y'^2+\mu _2y^{k_2}y^{\prime }=0
\end{gather}
is linearizable by
a generalized Sundman transformation.

For f\/inding the functions $F$ and $G$ one has to solve equations
(\ref{equa:07}), (\ref{equa:AAA})  and (\ref{equa:ZZZ}), which
become
\begin{gather*}
   F_{x} =0, \qquad F_{yy}=F_y(\mu_3y^{k_3+1}+k_2)/y, ~ G_x =0, \qquad G_y=
Gk_2/y.
\end{gather*}

For example, if $k_2=k_3$, one takes the simplest solution,
$F=\frac{1}{\mu_3}e^\frac{\mu_3y^{k_2+1}}{k_2+1}$ and $G=y^{k_2}$, and
the generalized Sundman transformation becomes
\begin{gather}
\label{ex:ZA} u =\frac{1}{\mu_3}e^\frac{\mu_3y^{k_2+1}}{k_2+1} ,\qquad dt = y^{k_2}dx.
\end{gather}
Equations (\ref{equa:08}), (\ref{equa:09}) and (\ref{equa:MN})
give
\begin{gather*}
    \beta =\mu_2,\qquad  \gamma=0,\qquad  \alpha=0.
\end{gather*}
Hence equation (\ref{ex:VV}) is mapped by the transformation
(\ref{ex:ZA}) into the linear equation
 \begin{gather*}
u''+\mu_2u'=0.
 \end{gather*}

If $\mu_3=0$, then equation (\ref{ex:00}) is
\begin{gather}
\label{ex:000}
y^{\prime \prime }+\mu _2y^{k_2}y^{\prime }+\mu _1y^{k_1}=0,
\end{gather}
where $\mu _2\neq 0$. The
Lie criteria \cite{Lie} show that the nonlinear equation (\ref{ex:000}) is linearizable by a~point transformation if and only if $k_1=3$, $%
k_2=1$ and $\mu _1=(\mu _2/3)^2$. In the particular case, $k_1=3$, $k_2=1$, $\mu
_1=1$ and $\mu _2=3$, one has the equation
\begin{gather}
\label{ex:aa}y^{\prime \prime }+3yy^{\prime }+y^3=0.
\end{gather}
Equation (\ref{ex:aa}) arises in many areas. Some of these are the analysis
of the fusion of pellets, the theory of univalent functions, the stability
of gaseous spheres, operator Yang--Baxter equations, motion of a free
particle in a space of constant curvature, the stationary reduction of the
second member of the Burgers hierarchy \cite{KarasuLeach}.
\begin{remark}
Equation (\ref{ex:aa}) is linearizable by a point
transformation and by a generalized Sundman transformation into
the equation $u^{\prime \prime }=0$ and $u^{\prime \prime
}+3u^{\prime }+2u=0$, respectively.
\end{remark}

Without loss of the generality\footnote{For example, scaling of the independent variable: $\overline{x}=\mu_2x$.},
one can assume that $\mu _2=1$. Hence, equation (\ref{ex:000}) becomes
\begin{gather}
\label{ex:bb}
y^{\prime \prime }+y^{k_2}y^{\prime }+\mu _1y^{k_1}=0.
\end{gather}
For this equation the coef\/f\/icients are
\begin{gather*}
\lambda _0=\mu _1y^{k_1},\qquad \lambda _1=y^{k_2},\qquad \lambda _2=0,\qquad \lambda
_3=k_2y^{k_2-1}, \\
\lambda _4=\mu _1k_1(k_1-1)y^{k_1-2}-k_2y^{2k_2-1},\qquad \lambda
_5=k_2y^{2k_2+1}.
\end{gather*}

If $k_2=0$, then $\lambda _5=0$ and equation (\ref{ex:bb}) is linearizable
by a generalized Sundman transformation.

If $k_2\neq 0$, then $\lambda _5\neq 0$ and conditions (\ref{equa:1b}), (\ref{equa:1a}), (\ref{equa:1d}), (\ref{equa:1g}) are reduced to
\begin{gather}
\label{eq:aug0509.1}\mu _1(2k_2+1-k_1)(k_2-k_1)=0.
\end{gather}
If conditions (\ref{eq:aug0509.1}) are satisf\/ied, then equation (\ref{ex:bb}) is
linearizable by a generalized Sundman transformation. Notice that in
the case $\mu _1(k_2-k_1)=0$, equation (\ref{ex:bb}) is trivially integrated
by using the substitution $y^{\prime }=H(y)$. A nontrivial case is $k_1=2k_2+1$.
In this case the functions~$F$ and~$G$ are solutions of the compatible overdetermined
system of equations
\begin{gather}
\label{eq:aug0509.2}F_x=0,\qquad F_{yy}=k_2F_y/y,\qquad G_x=0,\qquad G_y=k_2G/y.
\end{gather}
The general solution of equations (\ref{eq:aug0509.2}) depends on the value of the
constant $k_2$. For example, if $k_2\neq -1$, then a particular solution of system
(\ref{eq:aug0509.2}) is
\[
F=y^{k_2+1},\qquad G=y^{k_2}.
\]
Thus, the generalized Sundman transformation reduces equation (\ref
{ex:bb}) into the linear equation
\begin{gather*}
u^{\prime \prime }+u^{\prime }+(\mu _1(k_2+1))u=0.
\end{gather*}
\end{example}

\begin{remark}
Since equations (\ref{ex:01}), (\ref{ex:00}) and (\ref{ex:000})
are autonomous, their order can be reduced by the substitution $y'=f(y)$. It is worth to note that for
equations (\ref{ex:00}) and (\ref{ex:000})
the dif\/f\/iculties in using the generalized Sundman transformation are similar to solving the original equation by this reduction.
 \end{remark}

\section{Conclusion}

Application of the generalized Sundman transformation for the
linearization problem was analyzed in the paper. Since the method is well-known,
the ef\/f\/iciency of the method is not discussed in the paper.
The paper just warns that a researcher has
to be careful when using the well-known method for the linearization problem.
In particular, our examples show that, in
contrast to point transformations (S.~Lie results), for a
linearization problem via the generalized Sundman transformation
one needs to use the general form of a linear second-order
ordinary dif\/ferential equation instead of the Laguerre form.

\subsection*{Acknowledgements}
\label{Acknowledgements}

This research was supported by the Royal Golden Jubilee Ph.D.
Program of Thailand (TRF).

\pdfbookmark[1]{References}{ref}
\LastPageEnding

\end{document}